\documentstyle[amssymb,amsfonts]{amsart}

\newcommand\R{{{\mathbf R}}}

\newcommand\eps{{\varepsilon}}
\newenvironment{proof}{\noindent {\bf Proof} }{\endprf\par}
\def \endprf{\hfill  {\vrule height6pt width6pt depth0pt}\medskip}
\def\emph#1{{\it #1}}
\def\textbf#1{{\bf #1}}

\theoremstyle{plain}
  \newtheorem{theorem}[subsection]{Theorem}
  
  \newtheorem{proposition}[subsection]{Proposition}
  
  \newtheorem{corollary}[subsection]{Corollary}

\theoremstyle{remark}

\theoremstyle{definition}

\include{psfig}

\begin{document}

\title[Global regularity for logarithmically supercritical NLW]{Global regularity for a logarithmically supercritical defocusing nonlinear wave equation for spherically symmetric data}
\author{Terence Tao}
\address{Department of Mathematics, UCLA, Los Angeles CA 90095-1555}
\email{ tao@@math.ucla.edu}

\vspace{-0.3in}
\begin{abstract}
We establish global regularity for the logarithmically energy-supercritical wave equation $\Box u = u^5 \log(2+u^2)$ in three spatial dimensions for spherically symmetric initial data, by modifying an argument of Ginibre, Soffer and Velo \cite{gsv} for the energy-critical equation.  This example demonstrates that critical regularity arguments can penetrate very slightly into the supercritical regime.
\end{abstract}

\maketitle

\section{Introduction}

In this paper we consider the defocusing nonlinear wave equation
\begin{equation}\label{nlw}
\Box u = f(u)
\end{equation}
in three spatial dimensions, where $u: I \times \R^3 \to \R$ is a real scalar field (with time restricted to some time interval $I$), and $\Box u := - \partial_{tt} u + \Delta u$ is the d'Lambertian and\footnote{Our argument extends of course to other defocusing nonlinearities of $u^5 \log u$ type, but we select this one for sake of concreteness.} $f(u) := u^5 \log(2+u^2)$.  To avoid technicalities let us restrict attention to \emph{classical solutions} to this equation, by which we mean solutions which are infinitely smooth and are compactly supported in space for each fixed time $t$.  
Then one can easily verify that solutions to these equations have a conserved energy
\begin{equation}\label{energy}
E(u) := \int_{\R^3} \frac{1}{2} |\partial_t u |^2 + \frac{1}{2} |\nabla_x u|^2 + F(u)\ dx
\end{equation}
where $F$ is the nonlinear potential
$$ F(u) := \int_0^u v^5 \log(1+v^2)\ dv \sim u^6 \log(1+u^2).$$
Here and in the sequel we use $X \lesssim Y$ to denote the estimate $X \leq CY$ for some absolute constant $C$, and $X \sim Y$ to denote the estimate $X \lesssim Y \lesssim X$.  

The key feature to note is that this equation is just barely \emph{energy-supercritical}.  Indeed, the nonlinear component
$\int_{\R^3} F(u)\ dx \sim \int_{\R^3} u^6 \log(2+u^2)\ dx$ of the energy just barely fails to be controlled by the linear
component, in contrast to the energy-critical equation 
\begin{equation}\label{nlw-crit}
\Box u = u^5
\end{equation}
or the energy sub-critical equation $\Box u = |u|^{p-1} u$ for $1 < p < 5$).  There is substantial evidence to suggest that supercritical equations are, in general, quite badly
behaved; for instance, focusing supercritical equations $\Box u = -|u|^{p-1} u$ for $p>5$ can blow up instantaneously from finite energy initial data (see e.g. \cite{sogge:wave}), while the defocusing counterpart $\Box u = +|u|^{p-1} u$ is
highly unstable in the energy class \cite{lebeau}, \cite{cct}.  However, it seems that equations which are merely
\emph{logarithmically} supercritical are still barely within the range of the critical wellposedness theory.  We will
illustrate this phenomenon with a model result:

\begin{theorem}\label{main} 
Let $u_0, u_1 \in C^\infty(\R^3)$ be any spherically symmetric smooth initial data.  Then there is a unique global smooth solution to \eqref{nlw} with initial position $u(0,x) = u_0(x)$ and initial velocity $\partial_t u(0,x) = u_1(x)$.
\end{theorem}

In fact our arguments will also show that this equation is global wellposed for radial initial data in $H^2(\R^3) \times H^1(\R^3)$ and admits a scattering theory in this class, by repeating the arguments from \cite{gsv}; we shall omit the details as these consequences are rather standard.  One can in fact lower the regularity to $H^{1+\eps} \times H^\eps$ for any $\eps > 0$.  For the critical equation \eqref{nlw-crit}, the counterpart to Theorem \ref{main} was established in \cite{struwe} (with a particularly simple proof given later in \cite{gsv}), and then extended to arbitrary smooth initial data in \cite{g_waveI}, \cite{grillakis.semilinear}; see e.g. \cite{shatah-struwe} for a survey of these results.

Note from finite speed of propagation (and the classical local wellposedness theory for smooth solutions, see e.g. \cite{sogge:wave}) that we may restrict attention to spherically symmetric classical solutions (i.e. ones which are both smooth and compactly supported at any given time).  In particular we may justify all integration by parts computations, and all Sobolev and Lebesgue norms are finite (at least on compact time intervals $I$).

Our argument shall in fact be a small modification of that in \cite{gsv}.  The key point in that paper is that in the spherically symmetric case, one can use the Morawetz inequality to obtain a useful \emph{a priori} spacetime bound for solutions to \eqref{nlw}, namely
$$ \int_I \int_{\R^3} |u(t,x)|^8\ dx dt \lesssim E^2$$
where $E$ is the energy.  This estimate also holds for the supercritical equation \eqref{nlw-crit}, and in fact we have the slightly stronger estimate
\begin{equation}\label{icrit}
 \int_I \int_{\R^3} |u(t,x)|^8 \log(2 + |u(t,x)|^2)\ dx dt \lesssim E^2.
\end{equation}
Indeed, if we let $G(u) := uf(u) - 2F(u)$, then the standard Morawetz inequality (see e.g. \cite{shatah-struwe}) gives
\begin{equation}\label{gux}
 \int_I \int_{\R^3} \frac{G(u)}{|x|}\ dx dt \lesssim E.
\end{equation}
On the other hand, if we write $f(u) = u g(u)$ then we have the identity
$$ G(u) = \int_0^u 2v [g(u) - g(v)]\ dv;$$
since $g(u) = u^4 \log(2+u^2)$ is even and increasing in $u$ for $u$ positive, we then observe that
$$ G(u) \sim |u|^6 \log(2 + |u|^2).$$
Applying the standard radial Sobolev inequality\footnote{This is the one place where we rely (crucially) on the assumption of spherical symmetry.}
$$ |u(t,x)| \lesssim \frac{1}{|x|^{1/2}} \| \nabla u(t) \|_{L^2_x(\R^3)} \lesssim E^{1/2} / |x|^{1/2};$$
inserting these bounds into \eqref{gux} we obtain \eqref{icrit} as claimed.

We shall use \eqref{icrit}
to subdivide the time interval $I$ into subintervals on which the integral of $|u|^8 \log(1+u^2)$ is small, and obtain
good control on the $H^2(\R^3) \times H^1(\R^3)$ norm (say) on each of these subintervals in turn.  The argument is somewhat reminiscent of the famous result of Beale-Kato-Majda \cite{bkm} on global existence of smooth solutions for the 2D Euler equation; there the logarithm arises from a failure of the endpoint Sobolev embedding $H^1(\R^2) \not \subset L^\infty(\R^2)$ rather than from a logarithmically supercritical nonlinearity, but the general structure of the argument seems similar.  For instance, as in \cite{bkm} our final bounds will also be double-exponential in the initial data norms.

It is likely that some version of Theorem \ref{main} extends to non-radial initial data and to other slightly supercritical
equations.  For instance, in view of the quantitative control recently established in \cite{tao} for non-radial solutions to \eqref{nlw-crit} (but with exponential bounds rather than polynomial), it seems likely that the methods there can extend to an equation such as $\Box u = u^5 [\log \log(10 + u^2)]^c$ for some small $c > 0$.  We will not pursue these matters.

The author thanks Patrick G\'erard for suggesting this problem.

\section{Notation and Strichartz estimate}

We use $L^q_t L^r_x$ to denote the spacetime norm
$$ \| u \|_{L^q_t L^r_x(\R \times \R^3)} := (\int_\R (\int_{\R^3} |u(t,x)|^r\ dx)^{q/r}\ dt)^{1/q},$$
with the usual modifications when $q$ or $r$ is equal to infinity, or when the domain $\R \times \R^3$
is replaced by a smaller region of spacetime such as $I \times \R^3$.  We also adapt this notation to the Sobolev
spaces $H^s_x(\R^3)$ in the obvious manner.

For classical spherically symmetric solutions $u: I \times \R^3 \to \R$ to $\Box u = F$, we recall the \emph{Strichartz estimate}
\begin{equation}\label{strichartz}
\begin{split}
\| u \|_{L^2_t L^\infty_x(I \times \R^3)} &+ 
\| \nabla_{t,x} u \|_{L^\infty_t L^2_x(I \times \R^3)}\\ 
&\lesssim \| \nabla_{t,x} u(t_0) \|_{L^2_x(\R^3)} + \| F \|_{L^1_t L^2_x(I \times \R^3)}
\end{split}
\end{equation}
for any $t_0 \in I$.  This endpoint $L^2_t L^\infty_x$ would normally be forbidden, but for spherically symmetric solutions it is available (essentially thanks to the fundamental solution and the Hardy-Littlewood maximal inequality); see \cite{kl-mac}.
Many other such estimates are available, but these are the only ones we shall need.

\section{Main argument}

Let $u: I \times \R^3 \to \R$ be a classical solution to \eqref{nlw} with energy $E$, and let $t_0$ be the lower endpoint of $I$.  We consider the quantities
\begin{align*}
A &:= \int_I \int_{\R^3} |u(t,x)|^8 \log(2 + |u(t,x)|^2)\ dx dt \\
B&:= \sum_{j=0}^1 \| \nabla^j_x u \|_{L^2_t L^\infty_x(I \times \R^3)} + \| \nabla_{t,x} \nabla_x^j u \|_{L^\infty_t L^2_x(I \times \R^3)} \\
D &:= \| \nabla_{t,x} u(t_0) \|_{H^1_x(\R^3)}
\end{align*}
and consider the relationships between these quantities with each other.

From differentiating \eqref{nlw} $j$ times for $j=0,1$ and using \eqref{strichartz} and H\"older, we have
\begin{align*}
B &\lesssim \| \nabla_{t,x} u(t_0) \|_{H^1_x(\R^3)} + \sum_{j=0}^1 \| \nabla^j_x f(u) \|_{L^1_t L^2_x(I \times \R^3)} \\
&\lesssim D + \sum_{j=0}^1 \| u^4 |\nabla^j_x u| \log(2+u^2) \|_{L^1_t L^2_x(I \times \R^3)} \\
&\lesssim D + \| u \log(2 + u^2)^{1/8} \|_{L^8_t L^8_x(I \times \R^3)}^4 \sum_{j=0}^1 \| \nabla^j_x u \|_{L^2_t L^\infty_x(I \times \R^3)} \| \log(2+u^2) \|_{L^\infty_t L^\infty_x(I \times \R^3)}^{1/2} \\
&\lesssim D + A^{1/2} B \log^{1/2}(2 + \|u\|_{L^\infty_t L^\infty_x(I \times \R^3)}^2).
\end{align*}
Applying the Sobolev embedding 
$$ \| u \|_{L^\infty_t L^\infty_x(I \times \R^3)}
\lesssim \sum_{j=0}^1 \| \nabla_{t,x} \nabla^j_x u \|_{L^\infty_t L^2_x(I \times \R^3)}$$
we conclude that
\begin{equation}\label{cbound} B \lesssim D + A^{1/2} B \log^{1/2}(2 + B^2).
\end{equation}

From \eqref{cbound} and a simple continuity argument (increasing $I$ continuously starting from $t_0$) we thus have

\begin{proposition}\label{propcont}  Suppose that $A$ obeys the smallness condition
$$A \leq \frac{\epsilon_0}{\log(2+D)}$$ 
for some sufficiently small absolute constant $\epsilon_0 > 0$.  Then we have $B \lesssim D$.  In particular, if $t_1$ is the upper endpoint of $I$, then we have
$$ \| \nabla_{t,x} u(t_1) \|_{H^1_x(\R^3)} \leq C_0 \| \nabla_{t,x} u(t_0) \|_{H^1_x(\R^3)}$$
for some absolute constant $C_0 > 0$. 
\end{proposition}

We can iterate this as follows.

\begin{corollary}  Let $A,B,D$ be as above (with no smallness assumption on $A$).  Then we have
$$ B \lesssim (2 + D)^{(2+D)^{O(A)}}.$$
\end{corollary}

\begin{proof}  Define $D_n := C_0^n D$ for $n=0,1,\ldots$.  Observe (from the integral test) 
that for any $N \geq 1$ we have the inequality
$$ \frac{\epsilon_0}{\log(2+D_0)} + \frac{\epsilon_0}{\log(2+D_1)} + \ldots + \frac{\epsilon_0}{\log(2+D_N)}
\gtrsim \log[2 + \frac{N}{\log(2+D)}].$$
Thus we may find an integer
$$ 1 \leq N \lesssim (2 + D)^{O(A))}$$
such that
$$ \frac{\epsilon_0}{\log(2+D_0)} + \frac{\epsilon_0}{\log(2+D_1)} + \ldots + \frac{\epsilon_0}{\log(2+D_N)}
> A.$$
This allows us to subdivide the time interval $I$ as
$$ I = [t_0,t_1] \cup [t_1,t_2] \cup \ldots \cup [t_{N-1},t_N]$$
such that
$$ \int_{t_n}^{t_{n+1}} \int_{\R^3} |u(t,x)|^8 \log(2 + |u(t,x)|^2)\ dx dt \leq \frac{\epsilon_0}{\log(2+D_n)}$$
for all $0 \leq n < N$.  An easy induction using Proposition \ref{propcont} then shows that
$$ \| \nabla_{t,x} u(t_n) \|_{H^1_x(\R^3)} \leq D_n$$
and
$$ \sum_{j=0}^1 \| \nabla^j_x u \|_{L^2_t L^\infty_x([t_n,t_{n+1}] \times \R^3)} + \| \nabla_{t,x} \nabla^j_x u \|_{L^\infty_t L^2_x([t_n,t_{n+1}] \times \R^3)} \lesssim D_n
$$
for all $0 \leq n < N$.  Adding these estimates together we obtain the claim.
\end{proof}

Combining this with \eqref{icrit}, we see that
$$ B \lesssim (2+D)^{ (2+D)^{O(E^2)} }.$$
In particular, if $u: [0,T] \times \R^3 \to \R$ is a classical solution to \eqref{nlw}, we see that
the quantity $\| \nabla_{t,x} u \|_{L^\infty_t H^1_x([0,T] \times \R^3)}$ is bounded by a quantity depending only on the energy and $H^2_x \times H^1_x$ norm of the initial data, and which is independent of $T$.  In particular, by Sobolev embedding the $L^\infty_{t,x}$ norm on $[0,T] \times \R^3$ is also bounded uniformly in $T$. Classical existence theory (see e.g. \cite{sogge:wave}) then gives global regularity as $t \to +\infty$.  The analogous claim for $t \to -\infty$ then follows from time reversal symmetry.  This completes the proof of Theorem \ref{main}.
\endprf


\begin{thebibliography}{10}

\bibitem{bkm}
J. T. Beale, T. Kato, A. Majda, \emph{Remarks on the breakdown of smooth solutions for the 3-D Euler equations}, Comm. Math. Phys. \textbf{94} (1984), 61--66.

\bibitem{cct}
M. Christ, J. Colliander, T. Tao, \emph{Ill-posedness for nonlinear Schrodinger and wave equations}, to appear, Annales Institut Henri Poincar\'e.

\bibitem{gsv}
J. Ginibre, A. Soffer, G. Velo, \emph{The global Cauchy problem for the critical nonlinear wave equation}, Jour. Func. Anal., \textbf{110} (1992), 96--130.

\bibitem{ginebre:summarywave}
J. Ginibre, G. Velo, \emph{Generalized Strichartz Inequalities for the
Wave Equation}, Jour. Func. Anal., \textbf{133} (1995), 50--68.

\bibitem{g_waveI}
M.  Grillakis, \emph{Regularity and asymptotic behaviour of the wave equation with a critical nonlinearity},  Ann. of Math. \textbf{132}  (1990), 485--509.

\bibitem{grillakis.semilinear}
M. Grillakis, \emph{Regularity for the wave equation with a critical
nonlinearity}, Commun. Pure Appl. Math.,
\textbf{45} (1992), 749--774.

\bibitem{kap.energy}
L. Kapitanski, \emph{Global and unique weak solutions of
nonlinear wave equations}
Math. Res. Letters, \textbf{1} (1994), 211--223.

\bibitem{Keel-Tao} M. Keel, T. Tao, \emph{Endpoint Strichartz estimates},
Amer. J. Math., 120 (1998), 955--980.

\bibitem{kl-mac}
S. Klainerman, M. Machedon, \emph{Space-time Estimates for
Null Forms and the Local Existence Theorem}, Comm. Pure Appl.
Math., \underline{46} (1993), 1221--1268.

\bibitem{lebeau}
G.~Lebeau,
{\em Optique non lin\'eaire et ondes sur critiques},
S\'eminaire: \'Equations aux D\'eriv\'ees Partielles, 1999--2000,
Exp. No. IV, 13 pp., S\'emin. \'Equ. D\'eriv. Partielles,
\'Ecole Polytech., Palaiseau, 2000.

\bibitem{shatah-struwe}
J. Shatah, M. Struwe, \emph{Geometric Wave Equations}, Courant Lecture Notes in Mathematics \textbf{2} (1998).

\bibitem{sogge:wave}
C.~D. Sogge, \emph{Lectures on Nonlinear Wave Equations}, Monographs in Analysis
 II, International Press, 1995.

\bibitem{struwe}
M. Struwe, \emph{Globally regular solutions to the $u^5$ Klein-Gordon
equation},  Ann. Scuola Norm. Sup. Pisa Cl. Sci. \textbf{15}  (1988), 495--513.

\bibitem{tao}
T. Tao, \emph{Spacetime bounds for the energy-critical nonlinear wave equation in three spatial dimensions}, preprint.

\end{thebibliography}
\end{document}